\documentclass[plain]{article}
\usepackage{amsfonts,amssymb,amsmath}
\usepackage{multirow}
\usepackage{morefloats}
\usepackage{etex}
\usepackage[dvipsnames]{xcolor}
\usepackage{psfrag}
\usepackage{bm}
\usepackage[normalem]{ulem}
\usepackage{graphicx}

\usepackage{amsthm}
\definecolor{dgreen}{rgb}{0,0.5,0}
\newtheorem{thm}{Theorem}

\newtheorem{lemma}[thm]{Lemma}
\usepackage[english]{babel}

\begin{document}

%\begin{frontmatter}

\title{A five field formulation for flow simulations in porous media with fractures and barriers via an optimization based domain decomposition method}

\author{Stefano Scial\`o \footnote{stefano.scialo@polito.it}\\{\footnotesize Dipartimento di Scienze Matematiche, Politecnico di Torino,
		Corso Duca degli Abruzzi 24,}\\{\footnotesize 10129 Torino, Italy. Member of INdAM research group GNCS.}}

\maketitle

\begin{abstract}
The present work deals with the numerical resolution of coupled 3D-2D problems arising from the simulation of fluid flow in fractured porous media modeled via the Discrete Fracture and Matrix (DFM) model. According to the DFM model, fractures are represented as planar interfaces immersed in a 3D porous matrix and can behave as preferential flow paths, in the case of conductive fractures, or can actually be a barrier for the flow, when, instead, the permeability in the normal-to-fracture direction is small compared to the permeability of the matrix. Consequently, the pressure solution in a DFM can be discontinuous across a barrier, as a result of the geometrical dimensional reduction operated on the fracture.

The present work is aimed at developing a numerical scheme suitable for the simulation of the flow in a DFM with fractures and barriers, using a mesh for the 3D matrix non conforming to the fractures and that is ready for domain decomposition. This is achieved starting from a PDE-constrained optimization method, currently available in literature only for conductive fractures in a DFM. First, a novel formulation of the optimization problem is defined to account for non permeable fractures. These are described by a filtration-like coupling at the interface with the surrounding porous matrix. Also the extended finite element method with discontinuous enrichment functions is used to reproduce the pressure solution in the matrix around a barrier.

The method is presented here in its simplest form, for clarity of exposition, i.e. considering the case of a single fracture in a 3D domain, also providing a proof of the well posedness of the resulting discrete problem. Four validation examples are proposed to show the viability and the effectiveness of the method.  
\end{abstract}

\subsection*{Keywords}
Darcy flows; non conforming FEM meshes; 2D-3D flow coupling; XFEM

%\MSC[2010] 65N30 \sep 65N50 \sep 68U20 \sep 86-08

%\end{frontmatter}

\newcommand{\D}{\mathcal{D}}
\newcommand{\Dc}{\mathring{\mathcal{D}}}
\newcommand{\n}{\bm{n}}
\newcommand{\NI}{\n_i}
\newcommand{\Dt}{\tilde{\mathcal{D}}}
\newcommand{\Hone}[1]{\mathrm{H}^1\mathrm{(}#1\mathrm{)}}
\newcommand{\Honezero}[1]{\mathrm{H}^1_0\mathrm{(}#1\mathrm{)}}
\newcommand{\Hrm}{\mathrm{H}}
\newcommand{\tr}[2]{\mathrm{tr}_{#1} \mathrm{(}#2\mathrm{)}}
\newcommand{\Fvec}{\bm{f}}
\newcommand{\K}{\mathbf{K}}
\newcommand{\jump}[1]{\left[\!\!\left[#1\right]\!\!\right]}

\newcommand{\Vrm}{\mathrm{V}}
\newcommand{\Lrm}{\mathrm{L}}
\newcommand{\Ucal}{\mathcal{U}}
\newcommand{\Hcal}{\mathcal{H}}
\newcommand{\Gmat}{{\bf{G}} }
\newcommand{\aD}{a^{\mathcal{D}}}
\newcommand{\ai}{a_i}
\newcommand{\aF}{a^F}
\newcommand{\TD}{\mathcal{T}^\D_{\delta_\D}}
\newcommand{\TFi}{\mathcal{T}^{i}_{\delta_{F_i}}}
\newcommand{\TGi}{\mathcal{T}^{i}_{\delta_{\Gamma_i}}}
\newcommand{\TS}{\mathcal{T}^{i}_{\delta_{S_m,i}}}
\newcommand{\ND}{\mathcal{N}_{h_\D}}
\newcommand{\NFi}{\mathcal{N}_{h_i}}
\newcommand{\NGi}{\mathcal{N}_{q_i}}
\newcommand{\NGt}{\mathcal{N}_{q}}
\newcommand{\NS}{\mathcal{N}_{u_i}^m}
\newcommand{\Jd}{\mathcal{J}}
\newcommand{\Cmat}{{\bf{C}}}
\newcommand{\Bmat}{{\bf{B}}}
\newcommand{\Nhf}{\mathcal{N}_{h_\Fcal}}
\newcommand{\Nht}{\mathcal{N}_h}
\newcommand{\NSt}{\mathcal{N}_u}
\newcommand{\NSip}{\mathcal{N}_{u_i}^+}
\newcommand{\NSi}{\mathcal{N}_{u_i}}
\newcommand{\Rmat}{{\bf{R}}}
\newcommand{\Amat}{{\bf{A}}}
\newcommand{\Dmat}{{\bf{D}}}
\newcommand{\Emat}{{\bf{E}}}
\newcommand{\Zmat}{{\bf{Z}}}
\newcommand{\Omat}{{\bf{O}}}
\newcommand{\Mmat}{{\bf{M}}}
\newcommand{\zero}{\bf{0}}
\newcommand{\Nw}{\mathcal{N}_w}
\newcommand{\Fcal}{\mathcal{F}}

%\begin{figure}
%\centering
%\includegraphics[width=0.7\textwidth]{Figure/split_facce_new.pdf}
%\caption{A simple DFM (left) and its front view (right) with nomenclature exemplification}
%\label{block_split}
%\end{figure}

\section{Introduction}
Flows in porous media are almost ubiquitous in the description of natural phenomena, appearing, for example both in geological and biological contexts, and this makes efficient and reliable simulations of great relevance for a large variety of applications. Fractures in a porous media can have a significant impact on simulation results, such that, often, their presence needs to be taken into account to correctly capture flow properties. 
In this perspective, Discrete Fracture and Matrix (DFM) model \cite{FFNAdler,MJR2005} can be used, in which fractures in a porous media are explicitly modeled. Alternatives are represented by homogenization techniques \cite{Qi2005}, dual-porosity models \cite{DPbook}, and embedded discrete fracture matrix models \cite{Moinfar2014,TENE2017,Odsaeter2019,MANZOLI2021,LOSAPIO2023}. To reduce the computational cost of simulations, in DFMs fractures are dimensionally reduced to 2D planar interfaces embedded in a 3D porous matrix, and averaged constitutive equations for the flow are derived for the fractures. Additional equations are introduced at the fracture-matrix interface to close the problem \cite{MJR2005}. Depending on the hydraulic characteristics of the fractures and of the surrounding medium, different situations can be encountered. Indeed, fractures with a permeability higher than that of the surrounding porous medium act as preferential paths for the flow, whereas, fractures with low permeability might represent barriers. The pressure solution in the matrix domain is expected to be continuous across permeable fractures, with discontinuous fluxes as a consequence of flux leakage in the fracture plane. Around a low permeable barrier, instead, the pressure might be discontinuous, due to the presence of the interface that behaves as a physical barrier in this case.

One of the major sources of complexity in DFM simulation is related to the generation of a mesh conforming to the fractures, since, in practical applications, a large number of fractures might be present, forming an intricate network of intersections \cite{BPSf,Odsaeter2019,Antonietti2019,Chernyshenko2020}. On the other hand domain size call for the use of domain decomposition strategies to tackle computational issues \cite{FS13,formaggia2014,FKS,BDS}.

A large variety of numerical schemes is available in the literature for problems in porous media with fractures. As the matrix solution can be discontinuous across fractures, the choice of a mesh conforming to such interfaces is still attractive, resulting in more conventional discretization schemes. Recent approaches in this direction are proposed in \cite{CHEN2023}, where a Discontinuous Galerkin (DG) discretization is used or in \cite{HYMAN2022} where Mimetic Finite Differences (MFD) are chosen instead. In some approaches, a partial non-conformity (see \cite{FKS} for this kind of nomenclature) is allowed between the mesh of the porous medium and that of the fractures, in the sense that even if fractures can not cross the elements of the 3D mesh, the 2D mesh elements on each fracture can be arbitrarily placed with respect to the faces of the 3D mesh elements on the same plane. An example is in \cite{Chave2018}, where a Hibrid High Order scheme is proposed on this kind of mesh. Polygonal/polyhedral meshes have also been suggested as a possibility to easily generate conforming or partially conforming meshes of complex domains, as for example in \cite{Antonietti2016}, where MFD are used; in \cite{Antonietti2019} in conjunction with a DG scheme; or in \cite{BBFPSV,Fumagalli2019,Coulet2019,Benedetto2020} where the Virtual Element Method (VEM) is adopted. Other choices are also based on Finite Volume schemes, such as the two or multi-point flux approximation as in \cite{Sandve2012,Faille2016} and gradient schemes, as in \cite{Brenner2016}. 
As an alternative, we find methods allowing fully non conforming meshes for the 3D matrix and the lower dimensional fractures, requiring, however, ad-hoc discretization strategies. In \cite{Koppel2019} Lagrange multipliers are used to couple the matrix-fracture problems on arbitrary meshes, but the approach is limited to conductive fractures, with a solution in the porous matrix that is continuous across fractures. The Extended Finite Element Method (XFEM) can be used to reproduce irregular solution on non-conforming meshes \cite{XfemReview,XFEM3D}, and its use in the context of porous media flows can be found, e.g., in \cite{FS13,BPSa,formaggia2014,Schwenck2015,Flemisch2016,DELPRA2017}.

Here an approach for 3D-2D coupled problems with discontinuous solutions on non conforming meshes is proposed. The method is an extension of the approach in \cite{BDS}, where DFM problems on non conforming meshes are solved via an optimization-based domain decomposition approach. The method in \cite{BDS} however, only considers continuous 3D solutions across the fractures, whereas here we allow for possibly discontinuous solutions. This is achieved by means of a novel five-field scheme for domain decomposition combined with the use of the XFEM. According to this approach three interface unknowns are added to the unknown pressure in the matrix and in the fractures, thus giving five independent fields. The interface variables allow to decouple the 3D problem from the 2D problem on each fracture and a global solution is then seen as the minimum of a cost functional, expressing the error in the matching at the interfaces. The method is presented here in the simplest case of a domain with a single fracture, to better highlight the novelties of the present approach and keeping the notation as compact as possible. The extension to the case of multiple intersecting fractures, resulting in the mixed dimensional problem described in \cite{Nordbotten2019,Boon2021} however, can be obtained by simply applying the ideas here proposed to the coupling on the lower dimensional domains. The purpose of the present work is to show the viability of the method and its well posedness. 

The manuscript is organized as follows: the problem considered in the present work is described in Section~\ref{Pdescr}, and its formulation as a five-field optimization based optimization problem is reported in Section~\ref{FFopt}. Section~\ref{FFoptDiscr} describes the discrete formulation of the proposed formulation, with a proof of its well posedness provided in Section~\ref{WellPosedness}.

\section{Problem description}
\label{Pdescr}

\begin{figure}
\centering
\includegraphics[width=0.55\textwidth]{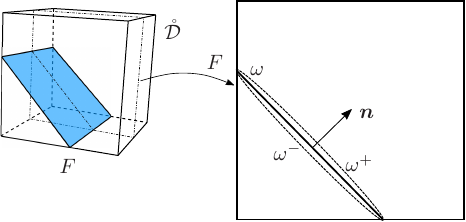}
\caption{Simple DFM example with nomenclature}
\label{dominio}
\end{figure}

Let us consider a 3D domain $\D$, representing a block of porous material crossed by a planar interface $F$, representing, instead a dimensionally reduced fracture in a DFM model. Let us further denote by $\Dc$ the domain $\D$ without the fracture, i.e. $\Dc:=\D\setminus F$. The boundary of $\D$, $\partial \D$, is split in a Dirichelet part $\Gamma_D$ and a Neumann part, $\Gamma_N$, such that $\Gamma_D\cup\Gamma_N=\partial\D$, $|\Gamma_D\cap \Gamma_N|=0$ and $|\Gamma_D|> 0$, whereas the boundary of $F$, $\partial F$, is subdivided in a Dirichelet part $\gamma_{D}$  and a Neumann part $\gamma_{N}=\partial F\setminus \gamma_D$. The boundary of $\Dc$, denoted by $\partial \Dc$, is instead split into three parts, $\partial \Dc:=\Gamma_D^\circ\cup \Gamma_N^\circ\cup \omega$, defined such that $\Gamma_D^\circ=\partial\Dc\cap\Gamma_D$, $\Gamma_N^\circ=\partial\Dc\cap\Gamma_N$ and $\omega=\partial\Dc\cap F$, this last part being the portion of the boundary of $\Dc$ that coincides with the fracture $F$. Denoting by $\n$ the unit normal to the fracture $F$, with fixed chosen orientation, we will further distinguish the two ``sides'' of $\omega$ as  $\omega^+$ and $\omega^-$, with $\omega^+$ being the portion of $\partial \Dc$ with outward pointing unit normal equal to $\n$. The nomenclature is exemplified in Figure~\ref{dominio}, where the three overlapping boundaries, $\omega$, $\omega^+$ and $\omega^-$ are shown separated for explanation purposes. 

The trace on a manifold $\Gamma$ of a sufficiently regular function $v$ is, in general, denoted as $\tr{\Gamma}{v}$. For brevity, however, the trace of a function $v$ on $\omega^\pm$ is denoted as $v_{\pm}$, for any sufficiently regular function $v$ defined on $\Dc$.

Following the approach in \cite{MJR2005}, the problem in primal formulation for the pressure distribution in $\D$ can be written as:\\
%\begin{minipage}{0.49\textwidth}
\begin{eqnarray}
-\nabla \cdot \left(K^\D \nabla H^\D\right)&=&g \quad \text{in} \ \Dc \label{MatrixEq} \\
-\nabla_F \cdot \left(K^F \nabla_F H^F\right) &=& -Q^+-Q^- \quad \text{in} \ F \label{FractureEq}\\
Q^\pm &=& -\eta (H^\D_{\pm}-H^F) \quad \text{on} \ \omega^\pm \label{Starling}\\
%\end{eqnarray}
%\end{minipage}
%\begin{minipage}{0.49\textwidth}
%\begin{eqnarray}
H^\D &=&0 \quad \text{on} \ \Gamma_D^\circ\\
K^\D \nabla H^\D \cdot \n_{\Gamma_N}&=&0 \quad \text{on} \ \Gamma_{N}^\circ \\
H^F &=&0 \quad \text{on} \ \gamma_{D}\\
K^F \nabla_F H^F \cdot \n_{\gamma_{N}}&=&0 \quad \text{on} \ \gamma_{N},  \label{Strong_neu}
\end{eqnarray}
%\end{minipage}\\
being $H^\D$ and $H^F$ the pressure distribution in $\Dc$ and in $F$, respectively. In the above equations $g$ is a source term defined in $\Dc$, $K^\D$ and $K^F$ are strictly positive definite tensors denoting, respectively, the (effective) hydraulic conductivity of the porous matrix and of the fracture on its tangential plane; $\eta$ is, instead, the effective hydraulic conductivity of the fracture in the direction normal to the fracture plane. The operator $\nabla$ is the three-dimensional gradient in $\Dc$ whereas $\nabla_F$ is the two-dimensional gradient on the plane containing the fracture. Homogeneous Dirichelet and Neumann boundary conditions are used for simplicity, with $\bm{n}_{\Gamma_N}$ the outward unit normal vector to $\Gamma_N$ and $\n_{\gamma_{N}}$ the outward unit normal vector to $\gamma_{N}$ on the plane of fracture $F$. 
Function $H^\D$ is allowed to have two different traces on $\omega^\pm$ and the source term in equation~\eqref{FractureEq} is related to the jump of the solution between the matrix and the fracture at the interface, as expressed by equation~\eqref{Starling}.

Let us introduce the function space $\Honezero{\Dc}=\left\lbrace v \in \Hone{\Dc}: \tr{\Gamma_D^\circ}{v}=0 \right \rbrace$, the space $\mathcal{Q}:=\Hrm^{\frac12}(\omega)$,  and the space $\Vrm=\left\lbrace v \in \Hone{F}, \tr{\gamma_{D}}{v}=0\right\rbrace$,
as well as bilinear forms $\aD : \Honezero{\Dc}\times \Honezero{\Dc} \mapsto \mathbb{R}$, $b:\mathcal{Q}\times\mathcal{Q} \mapsto \mathbb{R}$ and $\aF : \Vrm\times \Vrm \mapsto \mathbb{R}$, defined respectively as:
\begin{displaymath}
\aD\left(v, w \right)=\int_{\Dc} K^\D \nabla v \nabla w \, \mathrm{d}\Dc, \quad b(v,w)=\int_{\omega} \eta v w \, \mathrm{d}\omega,
\end{displaymath}
\begin{displaymath}
\aF\left(v, w \right)=\int_{F} K^F \nabla_F v \nabla_F w \, \mathrm{d}F.
\end{displaymath}

With the above definitions, the weak formulation of problem~\eqref{MatrixEq}-\eqref{Strong_neu} reads:
find $H^{\D} \in \Honezero{\Dc}$,  $H^F\in \Vrm$, such that:
\begin{eqnarray}
&&\aD\left(H^\D, v \right)+b\left(H^\D_{+},v_{+}\right)+b\left(H^\D_{-},v_{-}\right)= \nonumber \\
&& \hspace{2cm} \left(g,v\right)_{\Dc} +b\left(H^F,v_{+}\right)+b\left(H^F,v_{-}\right), \ \forall v\in \Honezero{\Dc} \label{Coupledweak1}\\
&&a^F\left(H^F,w\right)+2b\left(H^F,w\right)=b\left(H^\D_{+},w\right)\!+\!b\left(H^\D_{-},w\right)\!, \ \forall w\in \Vrm \label{Coupledweak2}
\end{eqnarray}
being $\left(v,w\right)_{d}$ the scalar product in $\mathrm{L}^2(d)$.

\section{Optimization based domain decomposition}
\label{FFopt}

The previous coupled system of equations can be de-coupled through the introduction of three additional interface variables, resulting in five independent unknowns for the problem. Equations~\eqref{Coupledweak1}-\eqref{Coupledweak2} are first re-written in terms of the new unknowns in a classical weak formulation, and subsequently as a PDE constrained optimization problem, which consists in the main novelty content of the present work.

Let us then introduce functions $\Psi^+$, $\Psi^- \in \mathcal{Q}$ representing the trace of $H^\D$ on $\omega^+$ and $\omega^-$, respectively, and function $\Psi^F\in\Vrm$ representing $H^F$. 
It is then possible to write problem~\eqref{Coupledweak1}-\eqref{Coupledweak2} as: find $H^{\D} \in \Honezero{\Dc}$, $H^F\in \Vrm$, $\Psi^+$,$\Psi^- \in \mathcal{Q}$ and $\Psi^F\in\Vrm$ such that, for all $v\in \Honezero{\Dc}$, for all $w\in \Vrm$:
\begin{eqnarray}
&&\aD\left(H^\D, v \right)+b\left(H^\D_+,v_{+}\right)+b\left(H^\D_-,v_{-}\right)= \nonumber \\
&& \hspace{2cm}\left(g,v\right)_{\Dc} +b\left(\Psi^F,v_{+}\right)+b\left(\Psi^F,v_{-}\right),  \label{Pweak1}\\
&&\aF\left(H^F,w\right)+ 2b\left(H^F,w\right)=b\left(\Psi^+,w\right)\!+\!b_i\left(\Psi^-,w\right)\!, \label{Pweak2}
\end{eqnarray}
with the following interface conditions:
\begin{eqnarray}
&\left\langle\Psi^+-H^\D_{+},\mu\right\rangle_{\mathcal{Q},\mathcal{Q}'}=0, & \forall \mu \in \mathcal{Q}', \label{Psi+} \\
& \left\langle\Psi^--H^\D_{-},\mu_i\right\rangle_{\mathcal{Q},\mathcal{Q}'}=0, & \forall \mu \in \mathcal{Q}',\label{Psi-} \\
& \left\langle\Psi^F-H^F,\kappa\right\rangle_{\Vrm,\Vrm'}=0, & \forall \kappa \in \Vrm'\label{PsiF}
\end{eqnarray}
where $\mathcal{Q}'$ and $\Vrm'$ denote the dual spaces of $\mathcal{Q}$ and $\Vrm$, respectively. It can be easily seen that problem~\eqref{Pweak1}-\eqref{PsiF} is equivalent to \eqref{Coupledweak1}-\eqref{Coupledweak2}, as conditions~\eqref{Psi+}-\eqref{PsiF} imply the (weak) equivalence of the newly introduced interface variables $\Psi^\pm$ and $\Psi^F$ with $H^\D_{\pm}$ and $H^F$, respectively. The advantage of the previous formulation consists in the possibility of fully de-coupling the problem on the fracture from the problem in the bulk domain. This is particularly advantageous in presence of a large number of fractures. In order to take full advantage of such decoupling a PDE-constrained domain decomposition approach is proposed for the problem, in view of its numerical resolution.
A cost functional is introduced to measure the error in the fulfilment of equations~\eqref{Psi+}-\eqref{PsiF}, the global solution of the problem being now obtained as the minimum of the functional constrained by the constitutive equations~\eqref{Pweak1}-\eqref{Pweak2}.
The cost functional $J(H^\D,H^F,\Psi^+,\Psi^-,\Psi^F)$, is defined as:
\begin{eqnarray}
&& J:=\|\Psi^+-H^\D_{+}\|^2_{\mathcal{Q}}+\|\Psi^--H^\D_{-}\|^2_{\mathcal{Q}}+\|\Psi^F-H^F\|^2_{\Vrm}. \label{Functional}
\end{eqnarray}
The solution to problem \eqref{Pweak1}-\eqref{PsiF} is then obtained as the minimum of functional $J$ constrained by the PDE equations on the 3D domain and on the fracture: 
\begin{equation}
\min J(H^\D,H^F,\Psi^+,\Psi^-,\Psi^F) \ \text{constrained by \eqref{Pweak1}-\eqref{Pweak2}}  \label{Pmin}
%\text{constrained by}& \text{\eqref{Pweak1}}&  \label{min1} \\
%\text{and by} &  \text{}& \label{min2}
\end{equation}

\section{Discrete problem}
\label{FFoptDiscr}
\newcommand{\mesh}[2]{\mathcal{T}^{#1}_{#2}}
\newcommand{\mpara}[2]{\delta^{#1}_{#2}}
\newcommand{\dof}[2]{\mathcal{N}^{#1}_{#2}}
\newcommand{\bfun}[2]{\varphi^{#1}_{#2}}
\newcommand{\NN}{\mathcal{N}}
\newcommand{\MM}{\mathcal{M}}
\newcommand{\xx}{{\bf{x}}}
\newcommand{\Gvec}{{\bf{g}}}

The optimization based domain decomposition formulation allows for a great flexibility in terms of meshing and discretization choices for the involved variables. In order to derive the discrete counterpart of problem \eqref{Pmin}, domain $\Dc$ is replaced by the whole $\D$. Lower-case letters are used for the variables of the discrete problems, which are $h^\D$, $h^F$, $\psi^\pm$ and $\psi^F$. A mesh, irrespective of the presence of the fracture, is built in $\D$ for $h^\D$ and, independently on $F$ for $h^F$. Further, independent meshes are introduced on $\omega$ for the three interface variables $\psi^\pm$ and $\psi^F$. The discontinuous behavior of the solution $h^\D$ across the interface $\omega$ on the non-conforming mesh is described using the eXtended Finite Element Method (XFEM), as discussed later in additional details.

Let us denote by $\mathcal{T}_{\delta^\D}^\D$ the tetrahedral mesh on $\D$, by $\mathcal{T}_{\delta^F}^F$ the triangular mesh on $F$, and by $\mathcal{T}_{\delta_\omega^\pm}^{\omega^\pm}$ and by $\mathcal{T}_{\delta_\omega^F}^{\omega^F}$ the polygonal meshes on the interface for $\psi^\pm$ and $\psi^F$, respectively. Quantities $\delta^\D$, $\delta^F$, $\delta_\omega^\pm$ and $\delta_\omega^F$ denote the mesh parameters of the corresponding meshes and indicate the maximum size of the elements in the mesh. The discrete variables are defined as:
\begin{displaymath}
h^\D=\sum_{k=1}^{\NN^\D} h^\D_k \phi_k(\xx), \quad h^F=\sum_{k=1}^{\NN^F} h^F_k \varphi_k(\xx_F), \quad \psi^\star=\sum_{k=1}^{\MM^\star} \psi^\star_k \vartheta^\star(\xx_F),
\end{displaymath} 
with $\star=\{+,-,F\}$, and $\xx_F$ a local coordinate system on $\omega$.
The discrete functional is denoted by $J_\delta$ and is obtained replacing the discrete variables in \eqref{Functional} and evaluating the norms of such finite dimensional quantities in $L^2(\omega)$, i.e.:
\begin{eqnarray}
&& J_\delta=\|\psi^+-h^\D_+\|^2_{\Lrm^2(\omega)}+\|\psi^--h^\D_-\|^2_{\Lrm^2(\omega)}+\|\psi^F-h^F\|^2_{\Lrm^2(\omega)}. \label{DiscrFunctional}
\end{eqnarray}
If we set $ h^T=\left[(h^\D)^T, (h^F)^T\right]^T$ and $\psi^T=\left[(\psi^+)^T, (\psi^-)^T, (\psi^F)^T\right]^T$, functional $J_\delta$ can be expressed in matrix form as:
\begin{equation}
J_\delta=
\begin{bmatrix}
h^T & \psi^T
\end{bmatrix}
\Gmat
\begin{bmatrix}
h \\
\psi
\end{bmatrix}, \quad
%\end{displaymath}
%where 
%\begin{displaymath}
\Gmat=
\begin{bmatrix}
\Gmat^{\D} & \zero & \Emat^+ & \Emat^- & \zero \\
\zero & \Gmat^F & \zero & \zero & \Emat^F \\
(\Emat^+)^T & \zero & \Gmat^{\psi^+} & \zero & \zero \\
(\Emat^-)^T & \zero & \zero & \Gmat^{\psi^-} & \zero \\
\zero & (\Emat^F)^T & \zero & \zero & \Gmat^{\psi^F}
\end{bmatrix}.
\label{Gmatrix}
\end{equation}
Matrices $\Gmat^{\D}\in\mathbb{R}^{\NN^\D\times\NN^D}$ and $\Gmat^{F}\in\mathbb{R}^{\NN^F\times\NN^F}$, $\Gmat^{\psi^\star}\in \mathbb{R}^{\MM^\star\times\MM^\star}$, $\star=\{+,-,F\}$ in $J_\delta$ are: 
\begin{displaymath}
(\Gmat^{\D})_{k,\ell}=\int_{\omega} \phi_{k|\omega^+}\phi_{\ell|\omega^+} \mathrm{d}\omega +  \int_{\omega} \phi_{k|\omega^-}\phi_{\ell|\omega^-}\mathrm{d}\omega;  \quad (\Gmat^{F})_{k,\ell}=\int_{F} \varphi_k \varphi_{\ell} \mathrm{d}F; 
\end{displaymath}
\begin{displaymath}
(\Gmat^{\psi^\star})_{k,\ell}=\int_{\omega} \vartheta^\star_k\vartheta^\star_\ell \mathrm{d}\omega; 
\end{displaymath}
whereas matrices $\Emat^\pm\in\mathbb{R}^{\NN^\D\times \MM^\pm}$ and $\Emat^F\in\mathbb{R}^{\NN^F\times \MM^F}$ are set as
\begin{displaymath}
(\Emat^\pm)_{k,\ell}=-\int_{\omega} \eta \phi_{k|\omega^\pm} \vartheta^\pm_\ell \mathrm{d}\omega; \quad (\Emat^F)_{k,\ell}=-\int_{\omega}\eta \varphi_{k} \vartheta^F_\ell\mathrm{d}\omega.
\end{displaymath}
The matrix form of the discrete constraint equations can be obtained after defining $\Amat^{\D}\in\mathbb{R}^{\NN^D\times\NN^D}$ as:
\begin{displaymath}
(\Amat^{\D})_{k\ell}=\int_\D K^\D \nabla \phi_k \nabla \phi_\ell + \int_{\omega} \eta \phi_{k|\omega^+} \phi_{\ell|\omega^+}\mathrm{d}\omega + \int_{\omega} \eta\phi_{k|\omega^-} \phi_{\ell|\omega^-}\mathrm{d}\omega,
\end{displaymath}
matrix  $\Amat^F\in\mathbb{R}^{\NN^F\times\NN^F}$ as:
\begin{displaymath}
(\Amat^F)_{k\ell}=\int_{F} K^F \nabla_F \varphi_k \nabla_F \varphi_\ell \mathrm{d}F +2 \int_{F}  \eta \varphi_k \varphi_\ell \mathrm{d}F,
\end{displaymath}
and matrices $\Bmat^{\star}\in\mathbb{R}^{\NN^F\times\MM^\star}$, $\star={+,-}$ and  $\Bmat^{F}\in\mathbb{R}^{\NN^\D\times\MM^F}$, respectively, as:
\begin{displaymath}
(\Bmat^{\star})_{k\ell}=\int_{\omega}\eta \varphi_{k} \vartheta_\ell^\star \mathrm{d}\omega , \quad (\Bmat^F)_{k\ell}=\int_{\omega}\eta  \phi_{k|\omega^+}\vartheta^F_\ell \mathrm{d}\omega + \int_{\omega}\eta  \phi_{k|\omega^-}\vartheta^F_\ell \mathrm{d}\omega.
\end{displaymath}
The discrete form of problem~\eqref{Pmin} reads:
%\begin{eqnarray}
%\Amat^{h^\D} h^\D &=& \Fvec - \Bmat^{\psi^+} \psi^+ - \Bmat^{\psi^-} \psi^- + \left(\Bmat^{\psi^F}_++\Bmat^{\psi^F}_-\right)\psi^F \label{Pdiscr1} \\
%\Amat^{h} h &=& \Cmat u + \Dmat^{\psi^+}\psi^+ + \Dmat^{\psi^-}\psi^- \label{Pdiscr2}
%\end{eqnarray}
\begin{eqnarray}
&&\min J_\delta=\begin{bmatrix}
h^T & \psi^T
\end{bmatrix} \Gmat \begin{bmatrix}
h\\
\psi
\end{bmatrix} \label{Jdiscr} \quad \text{such that}\nonumber \\
&&\Amat^{\D} h^\D = \Gvec + \Bmat^{F}\psi^F \label{Pdiscr1} \\
&&\Amat^{F} h^F = \Bmat^{+}\psi^+ + \Bmat^{-}\psi^- \label{Pdiscr2} 
\end{eqnarray}
where $\Gvec$ is the array resulting from the source term.

The above optimization problem \eqref{Jdiscr}-\eqref{Pdiscr2} can be rewritten as an unconstrained optimization problem, by formally replacing the constraints in the functional. Setting:
\begin{displaymath}
\Amat=\begin{bmatrix}
\Amat^\D & \zero\\
\zero & \Amat^F
\end{bmatrix}, \quad
\Bmat= \begin{bmatrix}
\zero & \zero & \Bmat^F\\
\Bmat^+ & \Bmat^- & \zero 
\end{bmatrix}, \quad
\Gmat^h=
\begin{bmatrix}
\Gmat^\D & \zero\\
\zero & \Gmat^F
\end{bmatrix}
\end{displaymath}
\begin{displaymath}
\Gmat^\psi=
\begin{bmatrix}
\Gmat^{\psi^+} & \zero & \zero\\
\zero & \Gmat^{\psi^-} & \zero\\
\zero & \zero & \Gmat^{\psi^F}\\
\end{bmatrix}, \quad
\Emat= \begin{bmatrix}
\Emat^+ & \Emat^- & \zero \\
\zero & \zero & \Emat^F
\end{bmatrix}, \quad
\bf{b}=\begin{bmatrix} \zero \\ \Gvec \end{bmatrix}
\end{displaymath}
we have:
\begin{displaymath}
\Gmat=\begin{bmatrix}
\Gmat^h & \Emat\\
\Emat^T & \Gmat^\psi
\end{bmatrix}
\end{displaymath}
such that the expression of $J_\delta$ becomes
\begin{displaymath}
J_\delta=h^T \Gmat^h h +h^T\Emat\psi + \psi^T\Emat^T h + \psi^T \Gmat^\psi \psi.
\end{displaymath}
By formally replacing $h=\Amat^{-1}\left(\Bmat\psi+\bf{b}\right)$ in the previous expression we get: 
\begin{eqnarray*}
&& J_\delta=\psi^T\left(\Bmat^T\Amat^{-T}\Gmat^h\Amat^{-1}\Bmat + \Emat^T\Amat^{-1}\Bmat + \Bmat^T\Amat^{-T}\Emat +\Gmat^\psi\right)\psi \nonumber \\ 
&& \hspace{2cm}+ 2\psi^T\left(\Bmat^T+\Amat^{-T}\Gmat^h\Amat^{-1}+\Emat^T\Amat^{-1}\right){\Gvec} + {\bf{b}}^{T}\Amat^{-T}\Gmat^\chi\Amat^{-1}{\Gvec},
\end{eqnarray*}
and the gradient of $J_\delta$ with respect to $\psi$ at $\bar{\psi}$ is given by:
\begin{displaymath}
\nabla_\psi J_\delta=\Bmat^T \bar{\lambda} + \Emat^T\bar{h}+\Gmat^\psi\bar{\psi}, \quad \bar{h}=\Amat^{-1}\left(\Bmat\bar{\psi}+{\bf{b}}\right), \quad
\bar{\lambda}=\Amat^{-T}\left(\Gmat^\chi \bar{h}+\Emat\bar{\psi}\right).
\end{displaymath}
This allows to solve the optimization problem via the conjugate gradient scheme, starting from an initial guess on the unknown $\psi$. Thanks to the block-diagonal structure of matrix $\Amat$, the gradient $\nabla_\psi J_\delta$ can be computed efficiently, decoupling the 3D problem from the 2D problem, and in a matrix-free approach.  

\section{Well posedness of the discrete problem}
\label{WellPosedness}
Well posedness of the discrete problem derives from the non-singularity of the system of optimality conditions of the constrained minimization problem \eqref{Jdiscr}-\eqref{Pdiscr2}. 

Let us set $\NN=\NN^\D+\NN^F$, $\MM=\MM^+ + \MM^- + \MM^F$ and $\NN^\text{tot}=\NN+\MM$. Further, let us denote by $\chi\in\mathbb{R}^{\NN^\text{tot}}$ the array of the unknowns, i.e. $\chi^T=\begin{bmatrix}h^T & \psi^T \end{bmatrix}^T$, by $\lambda\in\mathbb{R}^{\NN}$ the array of Lagrange multipliers, and by $\Cmat\in\mathbb{R}^{\NN\times\NN^\text{tot}}$ and $\bf{b}\in\mathbb{R}^{\NN}$ the quantities:
\begin{equation}
\Cmat=\begin{bmatrix}
\Amat^\D & \zero &  \zero &  \zero & -\Bmat^F\\
\zero & \Amat^F &  -\Bmat^+ & -\Bmat^- & \zero
\end{bmatrix}, \quad 
{\bf{b}}=\begin{bmatrix}
\zero\\
\Gvec
\end{bmatrix}.
\label{Cmatrix}
\end{equation}
Then the system of optimality conditions for \eqref{Jdiscr}-\eqref{Pdiscr2} can be written as:
\begin{equation}
\begin{bmatrix}
\Gmat  &\Cmat^T\\
\Cmat & \zero
\end{bmatrix}
\begin{bmatrix}
\chi\\
\lambda
\end{bmatrix}=
\begin{bmatrix}
\zero\\
{\bf{b}}
\end{bmatrix}.
\label{KKTsystem}
\end{equation}
Non singularity of the above system follows from Lemma~\ref{lemma1} below, using well known results of quadratic programming (see Theorem 16.2 in \cite{NWsecED}).
\begin{lemma}
\label{lemma1}
Let matrix $\Cmat\in\mathbb{R}^{\NN \times \NN^{\text{tot}}}$ be as in \eqref{Cmatrix},  and matrix $\Gmat\in\mathbb{\NN^{\text{tot}}\times\NN^{\text{tot}}}$ be as in \eqref{Gmatrix}, then $\Cmat$ is full row rank and $\ker{(\Gmat)}\cap\ker{(\Cmat)}=\zero$.
\end{lemma}
\begin{proof}
Non singularity of matrices $\Amat^\D$ and $\Amat^F$ derives from standard arguments, and thus it can be immediately shown that $\Cmat$ in \eqref{KKTsystem} is full row rank. Let us now take an element $\chi_0\in\ker{(\Gmat)}$. Then, by definition of $J_\delta$ we have that 
\begin{eqnarray*}
0&=&\chi_0^T\Gmat\chi_0=
\begin{bmatrix}
h_0^T & \psi_0^T
\end{bmatrix}
\Gmat
\begin{bmatrix}
h_0\\
\psi_0
\end{bmatrix} \nonumber \\
&=&\|\psi_0^+-\tr{\omega^+}{h_{0}^\D} \|^2_{\mathrm{L}^2(\omega)}+\|\psi_0^--\tr{\omega^-}{h_{0}^\D} \|^2_{\mathrm{L}^2(\omega)}+\|\psi_0^F-h_{0}^F \|^2_{\mathrm{L}^2(\omega)},
\end{eqnarray*}
and thus $\psi_0^+=\tr{\omega^+}{h_{0}^\D}$, $\psi_0^-=\tr{\omega^-}{h_{0}^\D}$ and $\psi_0^F=h_0^F$. 
Each element of $\ker{(\Cmat)}$ satisfies the constraint equations with $\Gvec=0$. Using now $\chi_0$ in the constraints gives for any $v\in \mathrm{H}^1_0(\Dc)$, and for all $w\in V$:
\begin{eqnarray*}
&&\int_{\Dc} K^\D \nabla h^\D_0 \nabla v \, \mathrm{d}\Dc+\int_{\omega} \eta \tr{\omega^+}{h^\D_0} v_+\, \mathrm{d}\omega+\int_{\omega} \eta \tr{\omega^-}{h^\D_0} v_- \,\mathrm{d}\omega = \nonumber \\
&& \hspace{7cm} \int_{\omega} \eta \psi_0^F v_+ \,\mathrm{d}\omega+\int_{\omega}\eta \psi_0^F v_- \,\mathrm{d}\omega,\\
&& \int_{F} K^F \nabla_F h_0^F \nabla_F w \, \mathrm{d}F + 2\int_F\eta h_0^F w \, \mathrm{d}F = \int_F\eta\left( \psi_0^+ w + \psi_0^-w\right) \, \mathrm{d}F.
\end{eqnarray*}
If we now take, in particular $v=h^\D_0$ and $w=h^F_0$, the above equations become:
%\begin{eqnarray*}
%&&\int_\D K^\D \nabla h^\D_0 \nabla h^\D_0 \, \mathrm{d}\D+\int_{\omega} \eta \psi_0^+ \psi_0^+\, \mathrm{d}\omega+\int_{\omega} \eta \psi_0^- \psi_0^- \,\mathrm{d}\omega =\int_{\omega} \eta \psi_0^F \psi_0^+ \,\mathrm{d}\omega+\int_{\omega}\eta \psi_0^F \psi_0^- \,\mathrm{d}\omega\\
%&& \int_{F} K^F \nabla_F h_0^F \nabla_F h_0^F + 2\eta \psi_0^F \psi_0^F \, \mathrm{d}F = \int_F\eta\left( \psi_0^+ \psi_0^F + \psi_0^-\psi_0^F\right) \, \mathrm{d}F
%\end{eqnarray*}
%and thus
\begin{eqnarray*}
&&\int_{\Dc} K^\D \nabla h^\D_0 \nabla h^\D_0 \, \mathrm{d}\Dc+\int_{\omega} \eta \left(\psi_0^+-\psi_0^F\right) \psi_0^+\, \mathrm{d}\omega+\int_{\omega} \eta\left( \psi_0^--\psi_0^F\right) \psi_0^- \,\mathrm{d}\omega =0\\
&& \int_{F} K^F \nabla_F h_0^F \nabla_F h_0^F \mathrm{d}F+\int_F\eta\left(\left(\psi_0^F-\psi_0^+\right) +\left(\psi_0^F-\psi_0^-\right)\right)\psi_0^F \, \mathrm{d}F = 0,
\end{eqnarray*}
and, by summing we get, ($F\equiv\omega$):
\begin{eqnarray*}
&&\int_{\Dc} K^\D \nabla h^\D_0 \nabla h^\D_0 \, \mathrm{d}\Dc +\int_{F} K^F \nabla_F h_0^F \nabla_F h_0^F \mathrm{d}F \\
&& + \int_{\omega} \eta \left(\psi_0^+-\psi_0^F\right) \left(\psi_0^+-\psi_0^F\right)\, \mathrm{d}\omega +\int_{\omega} \eta\left( \psi_0^--\psi_0^F \right) \left(\psi_0^--\psi_0^F\right) \,\mathrm{d}\omega =0\\
\end{eqnarray*}
and thus, for $\eta>0$, $\psi_0^+=\psi_0^F=\psi_0^-$ and consequently $h_0^\D=0$ and $h_0^F=0$.
\end{proof}

\section{Discretization choices}
\label{discrChoices}
The discretization of $h^\D$ is obtained by using linear Lagrangian finite elements on 3D tethrahedral meshes with additional basis functions used to describe the discontinuous behavior of the solution across the interface, following the XFEM paradigm. If a fracture entirely crosses the domain, then the Heaviside function $\mathcal{H}(\xx)$ can be used as enrichment, defined as: $\mathcal{H}(\xx)=\text{sign}\left(\n\cdot\left(\xx-\xx_0\right)\right)$, being $\n$ the unit normal to the fracture, and $\xx_0$ a fixed point on $F$. If, instead, a fracture terminates inside the domain, then the solution $h^\D$ is expected to be discontinuous across the fracture plane, but continuous around the fracture itself, and thus we need to choose an enrichment function with this kind of behaviour, i.e. discontinuous across the fracture but continuous elsewhere.
Let us consider here, for simplicity, a polygonal fracture $F$ with $n_e^\circ$ edges lying in the interior of the domain $\D$, consecutively numbered. Let then $\varepsilon_i(\xx_F)$ be the equation of the line containing the $i-$th edge in the fracture-local reference system, for $i=1,\ldots,n_e^\circ$, and $\mathcal{B}_F$ the affine mapping from the reference system in $\D$ to the fracture reference system, i.e. $\xx_F=\mathcal{B}_F(\xx)$. 
Then we can define function $\mathcal{E}(\xx)$ as:
\begin{displaymath}
\mathcal{E}(\xx)=
\begin{cases}
\sigma \prod_{i=1}^{n_e^\circ}\varepsilon_i^\varsigma(\mathcal{B}_F(\xx))& \text{if} \ \mathcal{B}_F(\xx)\in F\\
0 & \text{otherwise}
\end{cases},
\end{displaymath}
where $\sigma$ is a scaling parameter, chosen such that $\mathcal{E}(\xx)$ is equal to one in fracture barycenter, and the exponent $\varsigma$ can be chosen to control the shape of the enrichment function. In the proposed numerical example a value of $\varsigma=\frac12$ is used. We remark that $\mathcal{E}(\xx)$ is a bubble function on $F$ whenever all fracture edges are inside the domain, or the unit constant function when the fracture entirely crosses the domain, instead.  
We then choose function $\mathcal{H}(\xx)\mathcal{E}(\xx)$ as enrichment, observing that it is discontinuous across the fracture, but continuous (equal to zero) elsewhere in $\D$. This approach is similar to the one proposed in \cite{XFEM3D}, adapted to the present case.

\begin{figure}
\centering
\includegraphics[width=0.8\textwidth]{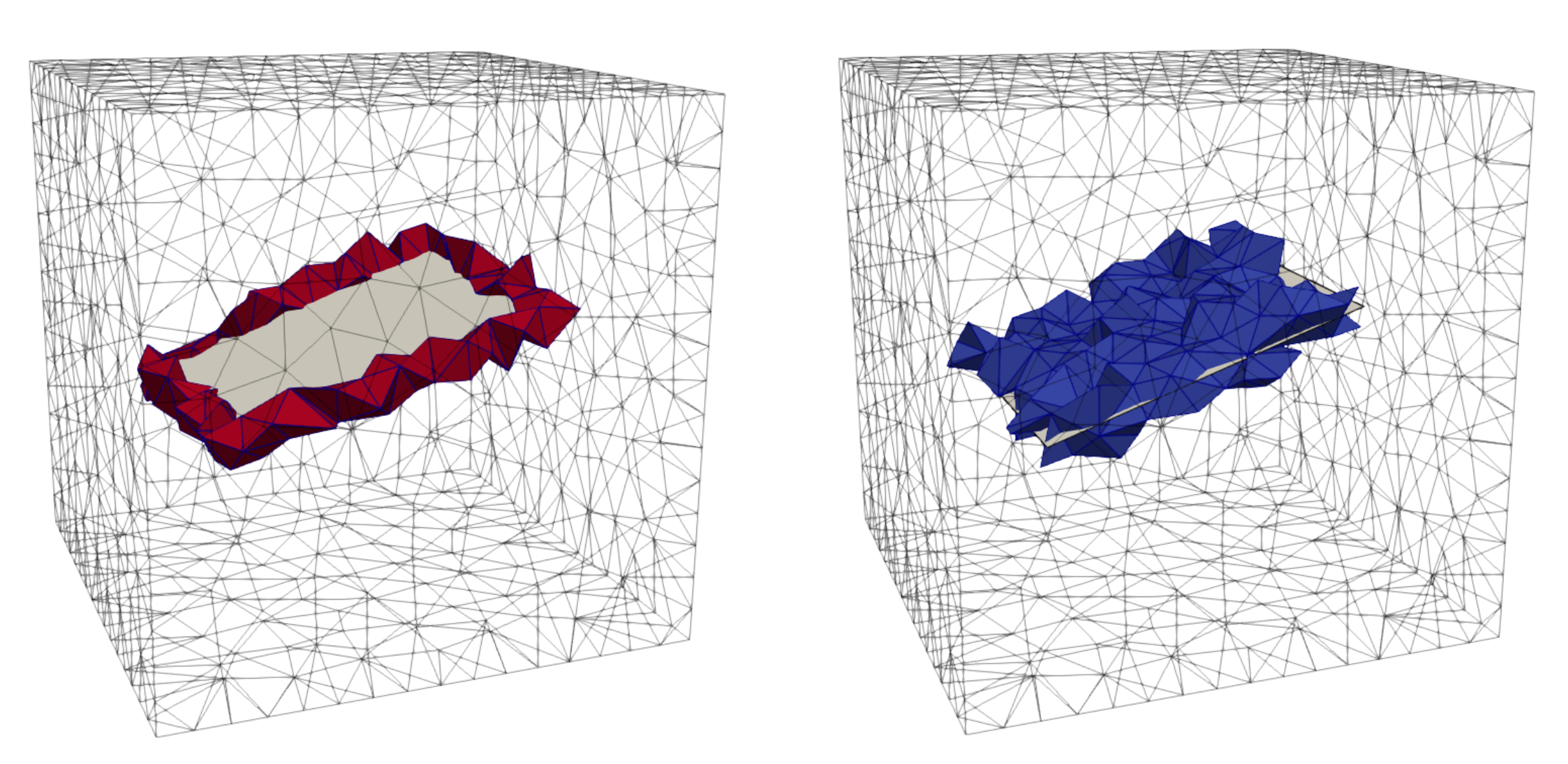}
\caption{Selection of elements to be enriched for an example embedded fracture: elements in  $\mathcal{T}^\D_\mathcal{E}$ in red on the left and elements in $\mathcal{T}^\D_\mathcal{H}$ in blue on the right.}
\label{enrEle}
\end{figure}

Let us denote by $\mathcal{T}^\D_\mathcal{H}$ the subset of the elements in $\mathcal{T}^{\D}_{\delta^\D}$ that are entirely cut by the fracture, and by $\mathcal{T}^\D_\mathcal{E}$ the set of mesh elements intersect by the edges of $F$ lying in the interior of $\D$, see Figure~\ref{enrEle}. Assuming, for simplicity of exposition, that the numbering of the vertexes is equal to that of the corresponding degrees of freedom (DOF), we can denote by $\mathcal{I}$ the set of all vertex indexes, and by $\mathcal{I}_\mathcal{E}$ the set of vertex indexes of elements in $\mathcal{T}^\D_\mathcal{E}$, and by $\mathcal{I}_\mathcal{H}$ the set of vertex indexes of elements in $\mathcal{T}^\D_\mathcal{H}$. Then the discrete function $h^\D$ can be written as:
\begin{eqnarray}
&& h^\D=\sum_{k\in\mathcal{I}} h^{s,\D}_k \phi_k^s(\xx) + \sum_{k\in\mathcal{I}_\mathcal{E}} h_k^{\mathcal{E},\D} \phi_k^{s}(\xx)(\mathcal{H}(\xx)\mathcal{E}(\xx)-\mathcal{H}(\xx_k)\mathcal{E}(\xx_k)) \nonumber \\ 
&& \hspace{4cm}+ \sum_{k\in\mathcal{I}_\mathcal{H}\setminus\mathcal{I}_\mathcal{E}} h_k^{\mathcal{I},\D} \phi_k^{s}(\xx)(\mathcal{H}(\xx)-\mathcal{H}(\xx_k));
\label{hdiscr}
\end{eqnarray}
where $\phi_k^s$ is the $k$-th linear Lagrangian basis function of standard FEM and $\xx_k$ are the coordinates of vertex $k$, with $\phi_\ell(\xx_k)=\delta_{\ell,k}$. The enrichment functions are shifted in order to be zero-valued in the nodes of the mesh. The interested reader is referred to the XFEM specific literature for further details.

%\begin{figure}
%\centering
%\includegraphics[width=0.8\textwidth]{Figure/EnrichmentFunction.pdf}
%\caption{Reconstruction of the global enrichment function for an example embedded fracture (zooming on the right).}
%\label{enrFun}
%\end{figure}

Standard linear Lagrangian finite element basis on the 2D triangular mesh $\mathcal{T}^F_{\delta^F}$ are used for the description of $h^F$. For simplicity, we choose, here, the same mesh for the fracture unknown and for the interface unknowns, i.e. $\mathcal{T}^{\omega^+}_{\delta^+}=\mathcal{T}^{\omega^-}_{\delta^-}=\mathcal{T}^{\omega^F}_{\delta^F}=\mathcal{T}^F_{\delta^F}$,  but we use piece-wise constant basis functions on such triangular mesh for $\psi^\star$, $\star=\{+,-,F\}$.

\section{Numerical results}
The present section is devoted to the presentation of the viability of the approach and its validation. The first two tests consider problems for which an analytic solution is known, and we measure the error between such solution and the solution obtained with the proposed approach. The last two examples, instead, give a qualitative comparison of the method with other reference solutions available in the literature. 

\begin{figure}
\centering
\includegraphics[width=0.9\textwidth]{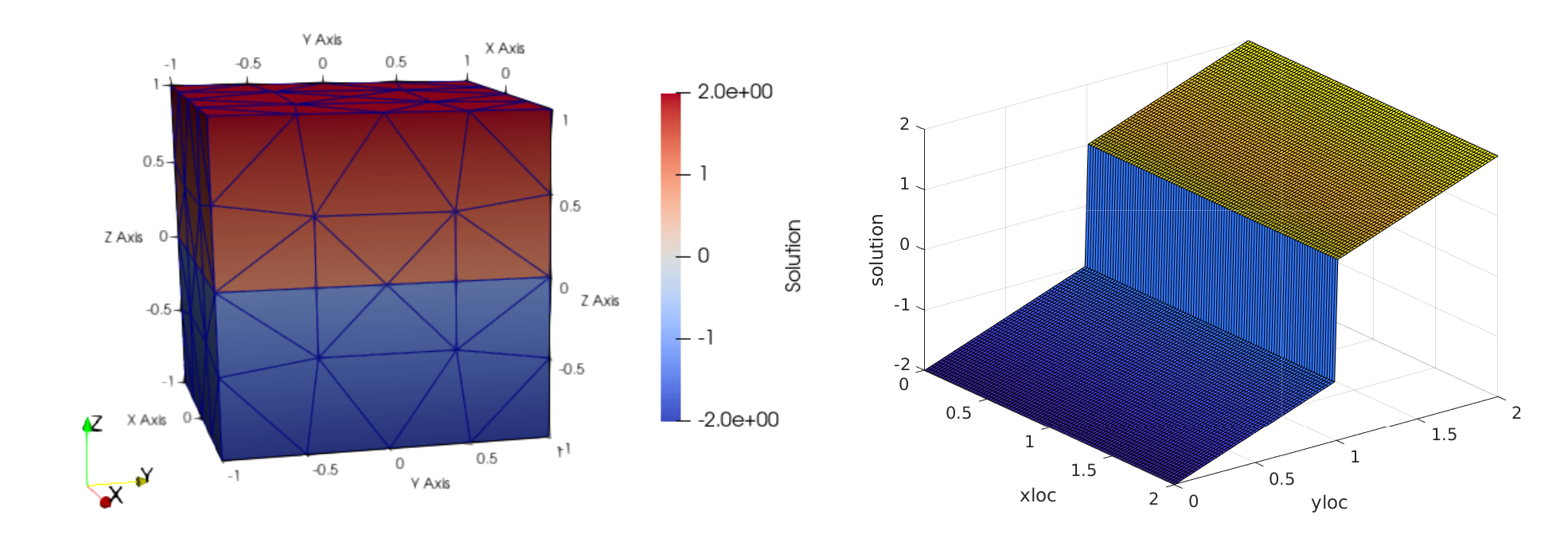}
\caption{\textit{test0}: solution on the 3D domain (left) and \textit{slice plot} of the solution on a plane orthogonal to $F$ (right).}
\label{test0}
\end{figure}

The first numerical experiment, named \textit{test0}, takes into account a problem with a solution that belongs to the discrete enriched FEM function space, as  described in Section~\ref{discrChoices}. The domain $\D$ is a cube with edge-length equal to $2$ and barycenter in the origin of a reference system $xyz$, see Figure~\ref{test0}, left. A fracture $F$ entirely crosses the domain, passing through the origin and having unit normal vector coinciding with the $z$-axis. Dirichelet boundary conditions equal to $-2$ are prescribed at the bottom face of the cube, i.e. on the face lying on the plane orthogonal to the $z$-axis and passing through point $P_0=[-1,-1,-1]$, and on its edges and vertexes. Dirichelet boundary conditions equal to $2$ are instead fixed on the top face, passing through point $P_1=[1,1,1]$ and on its edges and vertexes. Homogeneous Neumann boundary conditions are prescribed, instead, on the remaining part of the boundary of $\D$. Homogeneous Neumann boundary conditions are also applied on the whole boundary of $F$, and $K^\D=K^F=\eta=1$, $\Gvec=0$. The method is capable of obtaining the analytic solution
\begin{displaymath}
h^\D=
\begin{cases}
z+1 & z>0 \\
z-1 & z<0
\end{cases}, \quad h^F=0,
\end{displaymath}
independently of the meshsize. An example solution is reported in Figure~\ref{test0}, on the whole 3D domain (left) and on a plane orthogonal to the fracture (right). The 3D plot on the left represents the solution on sub-polyhedrons not crossing the interface, obtained cutting the original tetrahedral mesh along the fracture plane and evaluating the solution in the newly generated points (two different values on each side of the interface). All 3D plots in the remaining of the manuscript are obtained with the above procedure. The 2D plot on the right, instead, is obtained evaluating the discrete solution $h^\D(\xx)$ on $100\times100$ equally spaced points on the two orthogonal directions of a local reference frame on the plane. We will denote such plots as \textit{slice plots} in what follows. This kind of plots allows us to explicitly show the behavior of the solution, by reproducing the discrete function also inside mesh elements. It can be seen that, thanks to the presence of the enrichment basis function, the discrete solution is discontinuous across $F$, as expected. 

The second numerical test, labeled \textit{test1}, has the same domain geometry of \textit{test0}. Boundary conditions and the forcing term $\Gvec$ are chosen in such a way that the analytical solution is:
\begin{displaymath}
h^\D=
\begin{cases}
\text{e}^z & z>0\\
-\text{e}^{-z} & z<0\\
\end{cases}, \qquad
h^F=0;
\end{displaymath}
for $K^\D=K^F=\eta=1$. In particular we prescribe Dirichelet boundary conditions on the top and bottom faces, edges and vertexes of the domain, according to the chosen solution, and homogeneous Neumann boundary conditions on the remaining part of the boundary of $\D$ and on the whole $\partial F$. The numerical solution is computed on four meshes, with mesh parameter $\delta^\D$ ranging from $2\times10^{-2}$ to about $4\times10^{-5}$ and $\delta^F=\left(\delta^\D\right)^{\frac{2}{3}}$. The conjugate gradient (CG) scheme applied to the unconstrained optimization problem described at the end of Section~\ref{FFoptDiscr} is used to compute the numerical solution. Table~\ref{tab:test1} reports the number of degrees of freedom for $h^\D$ in column $\NN^\D$ and the DOFs for $h^F$  in column $\NN^F$ for the four chosen values of $\delta^\D$. The last column, $n_{\text{it}}$, shows the number of iterations of the CG algorithm to reach a relative residual of $10^{-7}$ for the various meshes: slightly more than $10$ iterations are sufficient in all cases.

\begin{table}
\centering
\begin{tabular}{cccc}
$\delta^\D$ & $\NN^\D$ & $\NN^F$ & $n_{\textit{it}}$\\
\hline
$2.0\times 10^{-2}$& 211 & 56 & 11 \\
$2.5\times 10^{-3}$& 1273 & 191 & 12 \\
$3.1\times 10^{-4}$& 8836 & 722 &  11\\
$3.9\times 10^{-5}$& 64754 & 2787 &  11\\
\end{tabular}
\caption{\textit{test1}: Number of DOFs for $h^\D$ ($\NN^\D$) number of DOFs for $h^F$ ($\NN^F$) and number of iteration of the CG algorithm ($n_{\text{it}}$) to reach a relative residual of $10^{-7}$, for different values of mesh parameter $\delta^\D$.}
\label{tab:test1}
\end{table}

\begin{figure}
\centering
\includegraphics[width=0.9\textwidth]{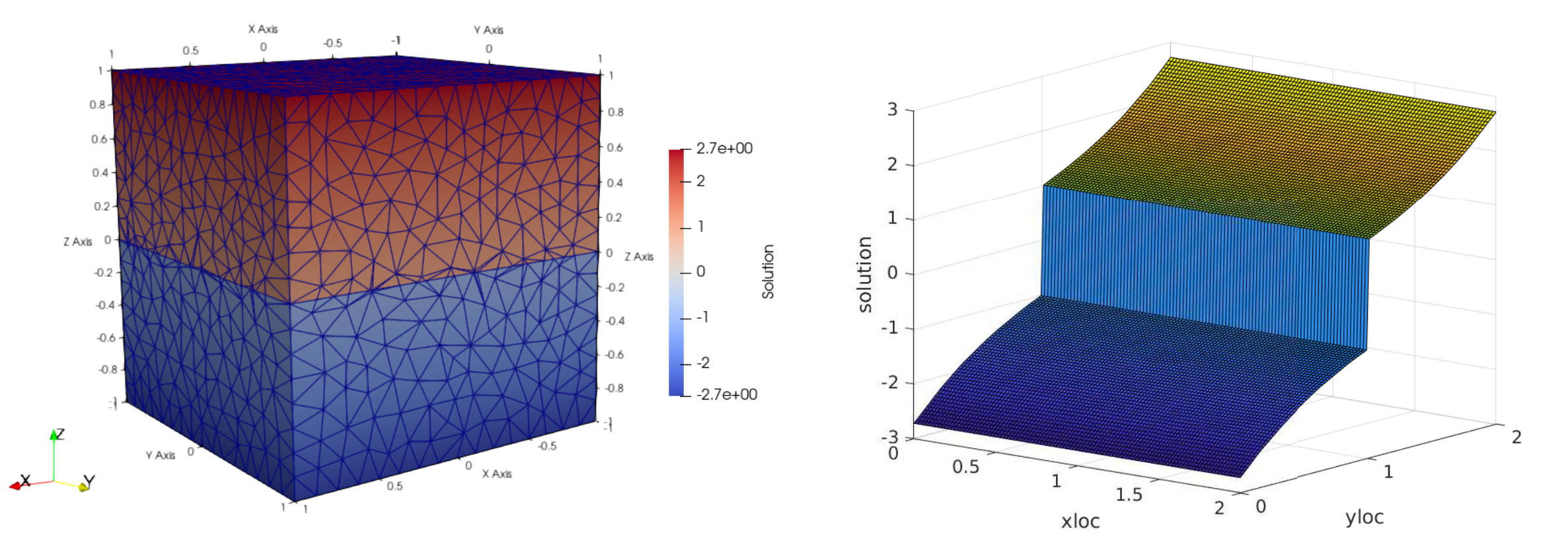}
\caption{\textit{test1}: solution on the 3D domain (left) and \textit{slice plot} of the solution on a plane orthogonal to $F$ (right)}
\label{fig:test1B}
\end{figure}

\begin{figure}
\begin{minipage}{0.49\textwidth}
\centering
\includegraphics[width=0.9\textwidth]{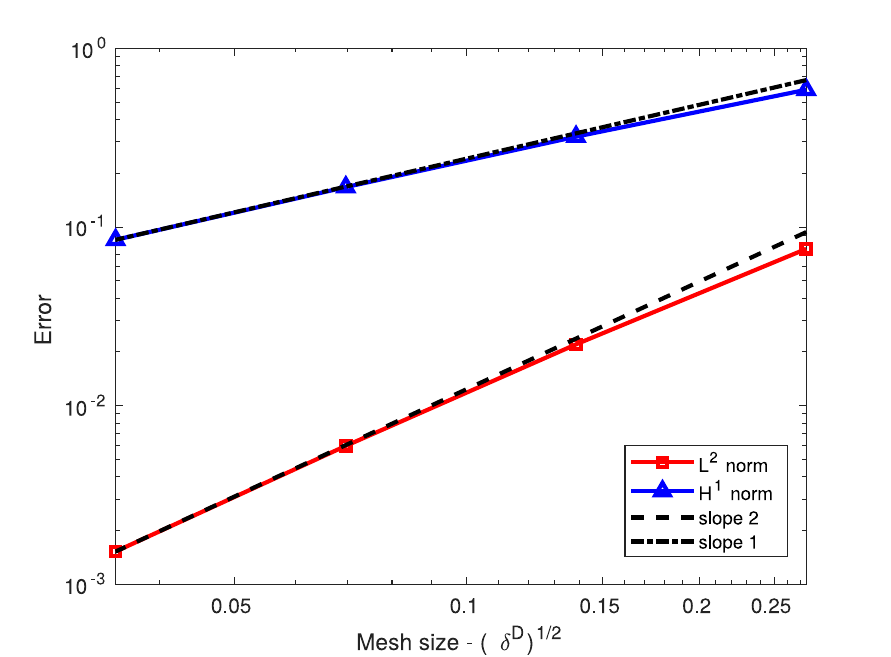}
\caption{\textit{test1}: Error convergence curves against mesh refinement}
\label{fig:test1A}
\end{minipage}
\begin{minipage}{0.49\textwidth}
\includegraphics[width=0.99\textwidth]{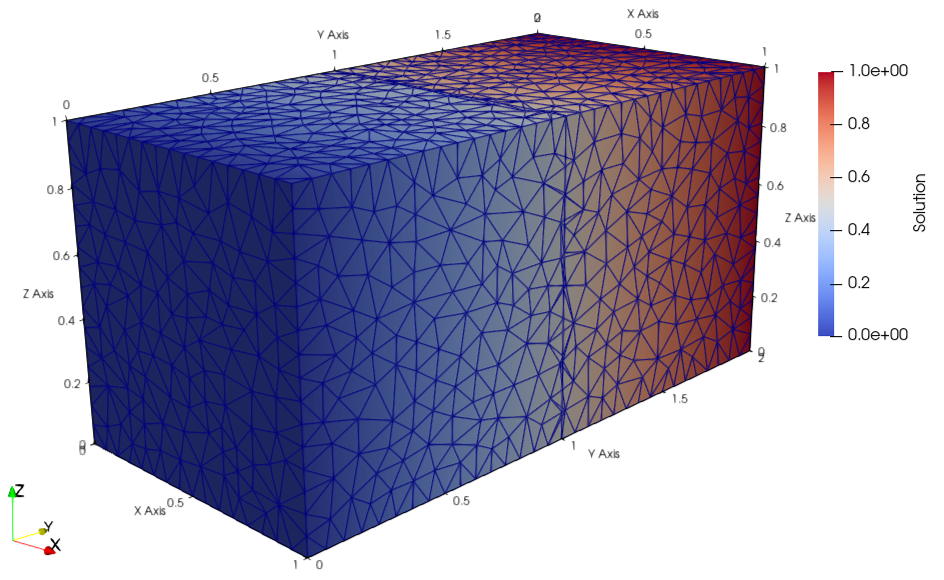}
\caption{\textit{test2}: solution on the 3D domain}
\label{fig:test2A}
\end{minipage}

\end{figure}

The solution obtained on the mesh with $\delta^\D=2.5\times 10^{-3}$ is displayed in Figure~\ref{fig:test1B}, left, on the whole domain , whereas its \textit{slice plot} on a plane orthogonal to the fracture can be found in Figure~\ref{fig:test1B}, right.
Figure~\ref{fig:test1A} displays convergence curves in the $L^2$ and $H^1$ norm of the error between the analytical and the numerical solution of the 3D problem. It can be noticed that the optimal trends for the chosen approximation space are obtained, thanks to the presence of the discontinuous basis functions, despite the non-conformity of the mesh. 

\begin{figure}
\centering
\includegraphics[width=0.99\textwidth]{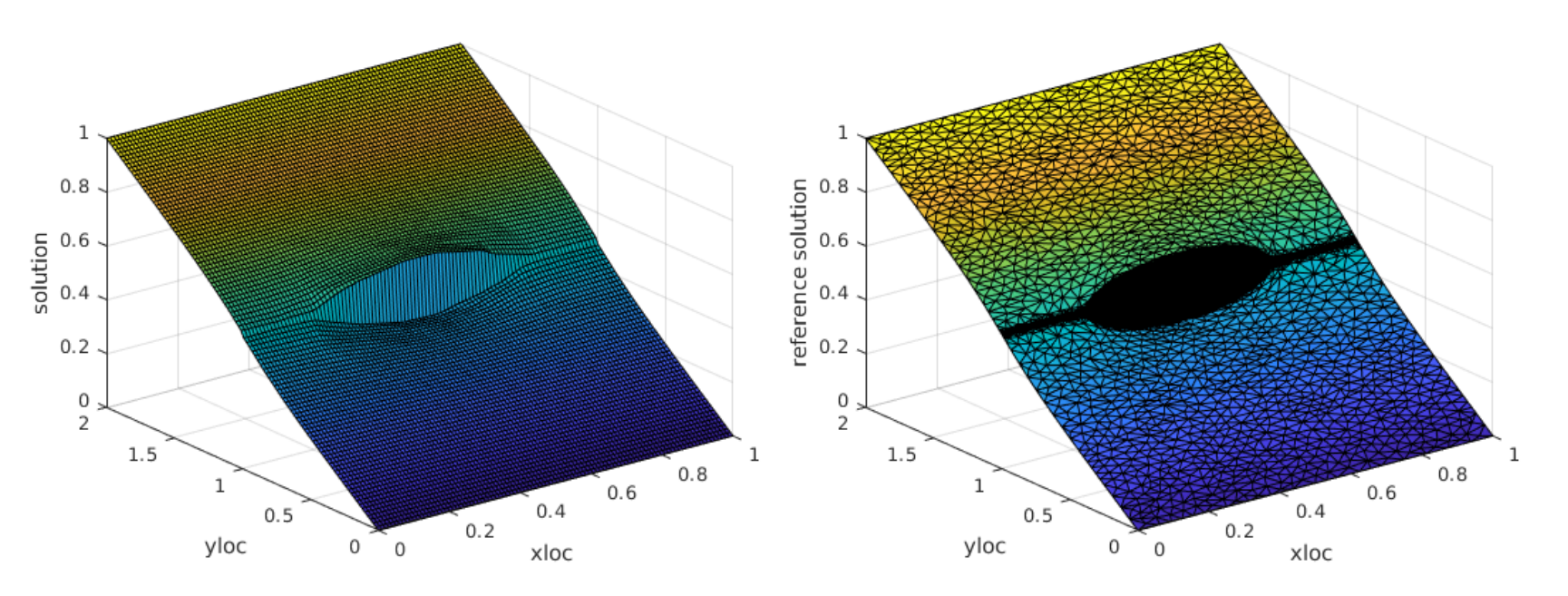}
\caption{\textit{test2}: \textit{Slice plot} of the solution with the proposed approach (left) and reference solution (right).}
\label{fig:test2B}
\end{figure}

\begin{figure}
\centering
\includegraphics[width=0.8\textwidth]{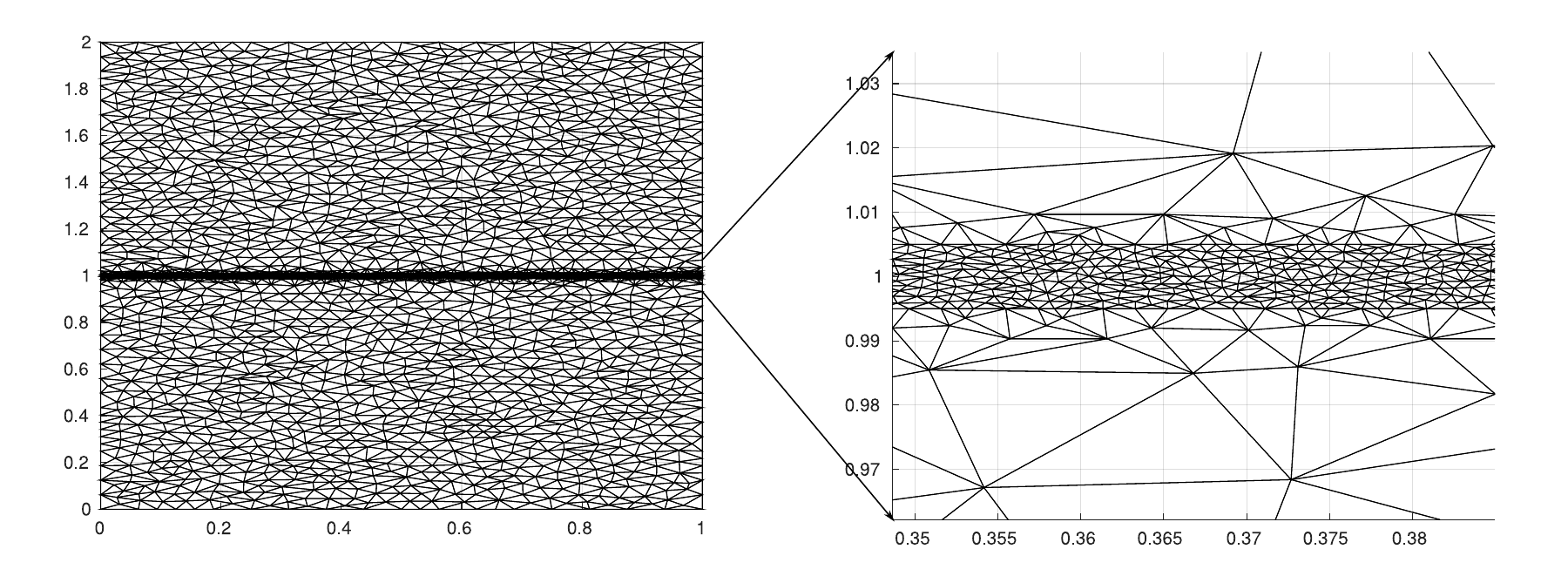}
\caption{\textit{test2}: Mesh for the reference solution with a zoom within fracture width}
\label{fig:test2BM}
\end{figure}

\begin{figure}
\begin{minipage}{0.45\textwidth}
\centering
\includegraphics[width=0.99\textwidth]{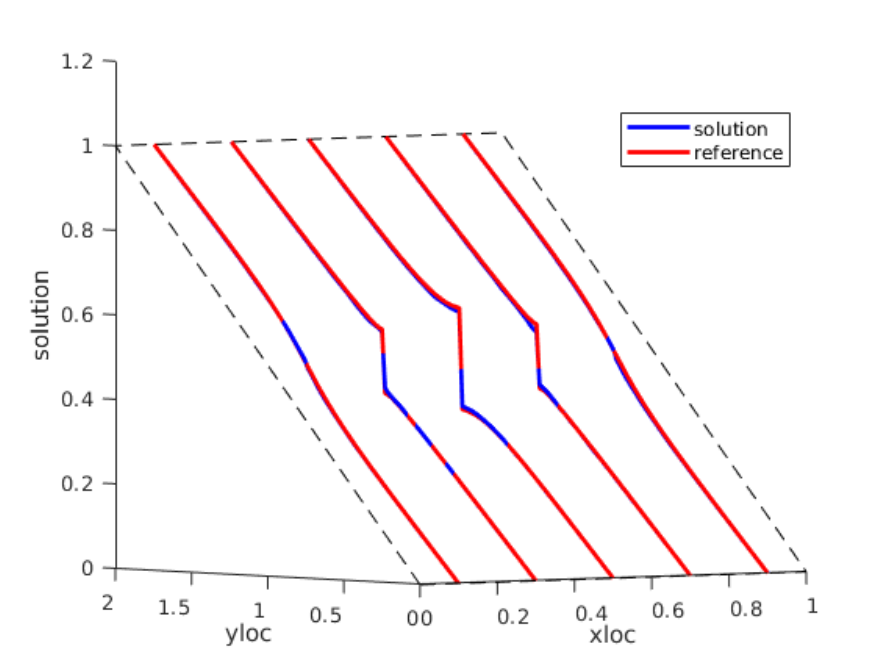}
\caption{\textit{test2}: Comparison of the reference solution and proposed solution on five segments on the $xy$ plane at $z=0.5$ for $x=\{0.1, 0.3,0.5, 0.7, 0.9\}$.}
\label{fig:test2C}
\end{minipage}
\begin{minipage}{0.04\textwidth}
~
\end{minipage}
\begin{minipage}{0.45\textwidth}
\centering
\includegraphics[width=0.99\textwidth]{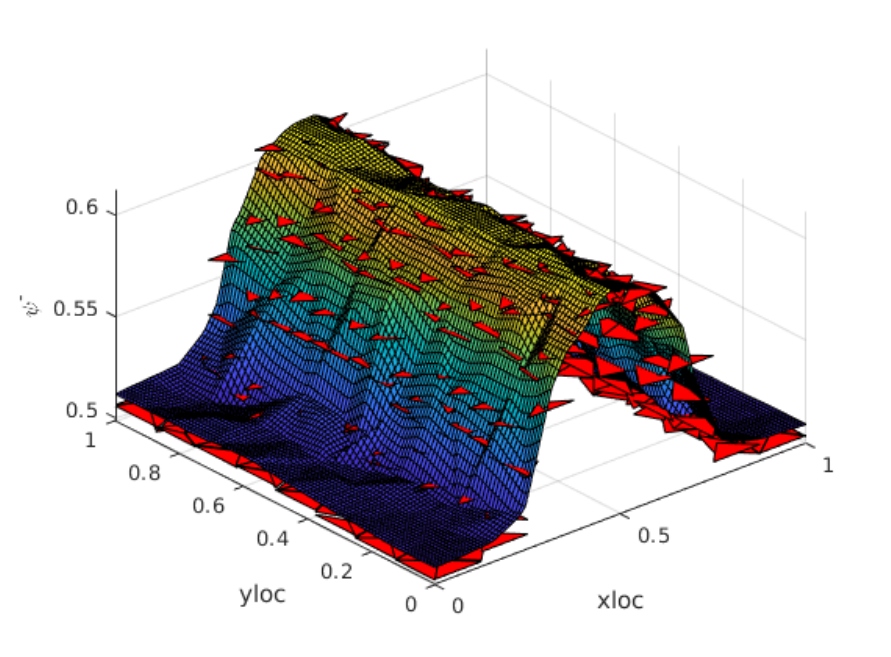}
\caption{\textit{test2}:\textit{Slice plot} of the solution on the plane $\omega^-$ compared to function $\psi^-$.}
\label{fig:test2D}
\end{minipage}
\end{figure}

The third numerical test, labelled \textit{test2}, is inspired by the second test case in Section 6.2 of \cite{MJR2005}. In the reference, a 2D problem is proposed, and here we consider a similar problem in three dimensions, obtained by uniformly extending in the third dimension the data available in \cite{MJR2005}.  The domain $\D$ is a cuboid with barycenter in $[0.5,1,0.5]$ and one vertex placed in the origin of a reference system $xyz$, see Figure~\ref{fig:test2A}. Edges parallel to the $y$-axis have length $2$, whereas all other edges have length $1$. Fracture $F$ lies on the $xz$ plane, passing through the barycentre, and entirely crossing the domain. A homogeneous Dirichelet boundary condition is set on cube face, edges and vertices lying on a plane orthogonal to the $y-$axis and passing through the origin, and a unitary Dirichelet boundary condition is instead set on face, edges and vertices lying on a plane orthogonal to the $y-$axis and passing through point $[1, 2, 1]$. Null Neumann boundary conditions are fixed on the remaining of $\partial\D$ and also on the whole boundary of $F$. In this example $K^\D=1$, whereas, being $d=10^{-2}$ the thickness of the original equi-dimensional fracture, it is:
\begin{displaymath}
K^F(\xx_F)=
\begin{cases}
d & x_F<0.25 \ \text{or} \ x_F>0.75\\
2\times 10^{-3}d & \text{otherwise}
\end{cases}
\end{displaymath}
and
\begin{displaymath}
\eta(\xx_F)=
\begin{cases}
1/d & x_F<0.25 \ \text{or} \ x_F>0.75\\
2\times 10^{-3}/d & \text{otherwise}
\end{cases}
\end{displaymath}
with $\xx_F=(x_F,y_F)$ a local reference system on $F$. As a consequence, fracture $F$ acts both as a permeable fracture where $K^F/d = K^\D$ and $\eta d= K^\D$, and as a barrier where, instead, $K^F/d$ and $\eta d$ are much smaller than $K^\D$.
A 3D mesh with about $2.5\times10^3$ degrees of freedom is considered, whereas the 2D meshes (the same for $h^F$ and all the interface variables) count $181$ DOFs for $h^F$ and $328$ DOFs for $\psi^\star$, $\star={+,-,F}$. $115$ CG iterations are required to reach a solution with a relative residual lower than $10^{-7}$. The obtained solution is reported in Figure~\ref{fig:test2A} in 3D. 

Let us denote by $F^\perp$ the rectangle given by the intersection of the plane $z=0.5$ with the domain $\D$ ($F^\perp$ is thus orthogonal to $F$). Figure~\ref{fig:test2B} shows, on the left, a \textit{slice plot} of the solution on $F^\perp$.  A 2D reference solution is also computed on $F^\perp$ for this problem. The reference solution is obtained using standard finite elements and explicitly representing the thickness of the fracture on $F^\perp$. The  triangular mesh for the reference solution is conforming to the fracture, is refined near the fracture, and inside fracture aperture has a maximum element diameter equal to $d/10$, as shown in Figure~\ref{fig:test2BM}. The reference mesh counts about $1.1\times10^4$ degrees of freedom, and the solution obtained on such mesh is shown in Figure~\ref{fig:test2B}, right, highlighting a good agreement with the solution obtained with the proposed method. In Figure~\ref{fig:test2C} a further qualitative comparison of the two solutions is proposed, through the overlapped plots of the solutions on five lines on $F^\perp$, placed at values of $x=\{0.1,0.3,0.5,0.7,0.9\}$, which again show a very good agreement. Finally, Figure~\ref{fig:test2D} proposes a comparison between the solution $h^-$ and the corresponding interface variable $\psi^-$. The matching between these two variables is given by the minimization of the functional, and we can notice a good agreement.

\begin{figure}
\centering
\includegraphics[width=0.5\textwidth]{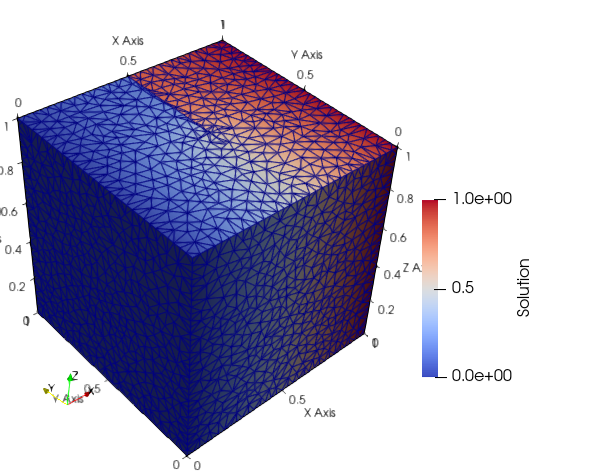}
\caption{\textit{test3}: solution on the 3D domain}
\label{fig:test3}
\end{figure}

\begin{figure}
\centering
\includegraphics[width=0.99\textwidth]{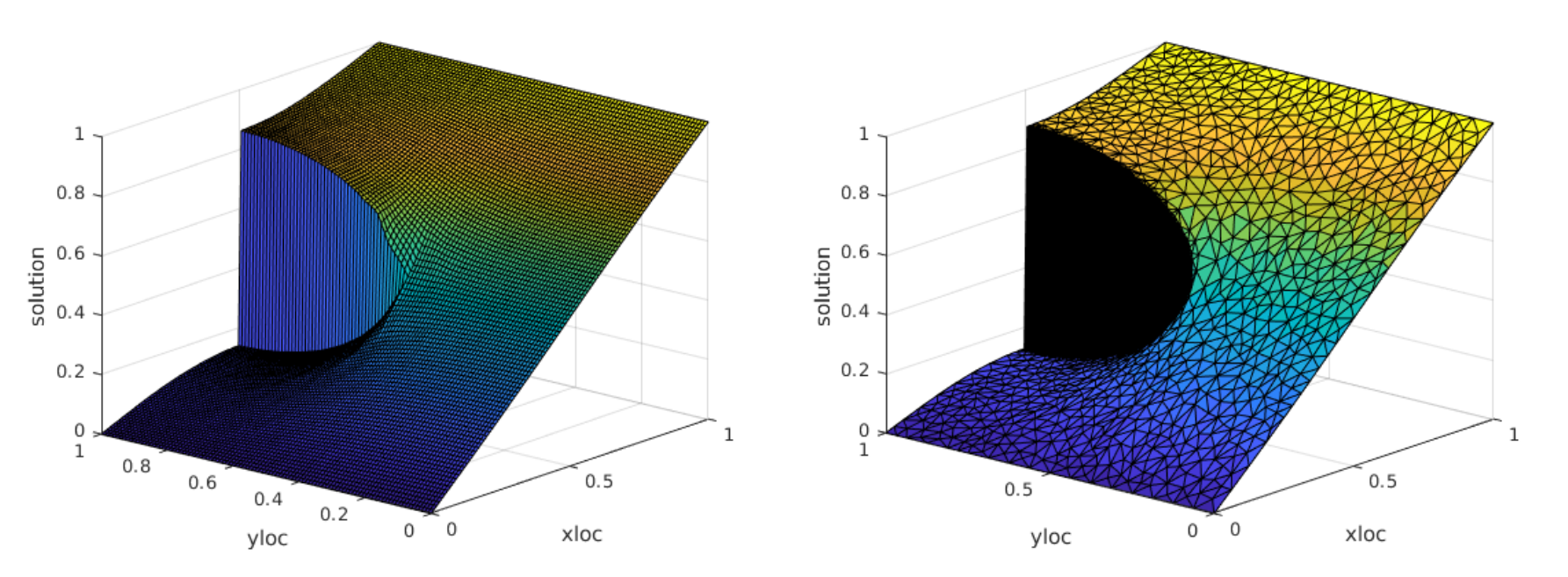}
\caption{\textit{test3}: \textit{Slice plot} of the solution with the proposed approach (left) and reference solution (right).}
\label{fig:test3A}
\end{figure}

\begin{figure}
\begin{minipage}{0.45\textwidth}
\centering
\includegraphics[width=0.99\textwidth]{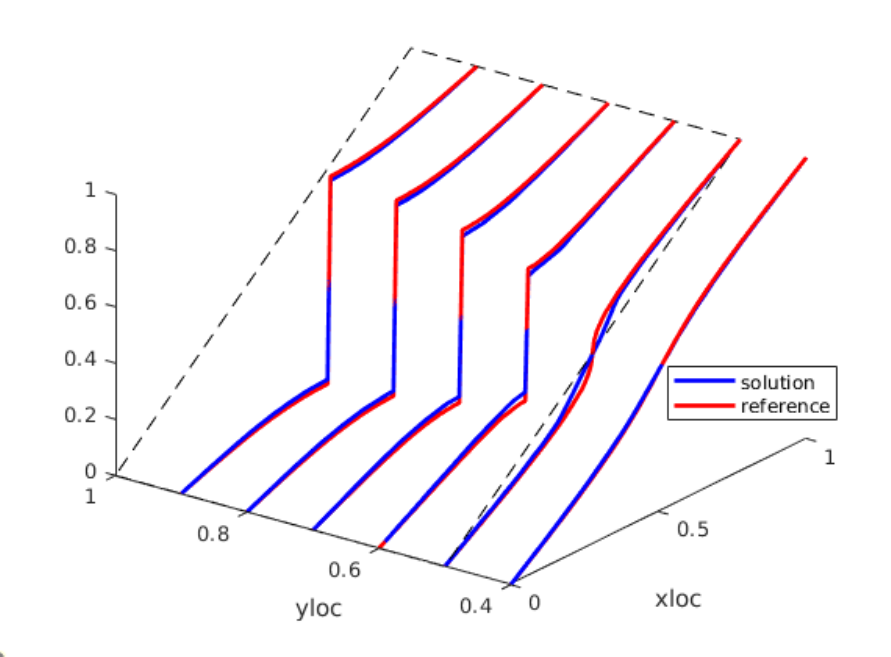}
\caption{\textit{test3}: Comparison of the reference solution and proposed solution on six segments on the $xy$ plane at $z=0.5$ for $y=\{0.4, 0.5,0.6, 0.7,0.8, 0.9\}$.}
\label{fig:test3B}
\end{minipage}
\begin{minipage}{0.04\textwidth}
~
\end{minipage}
\begin{minipage}{0.45\textwidth}
\centering
\includegraphics[width=0.99\textwidth]{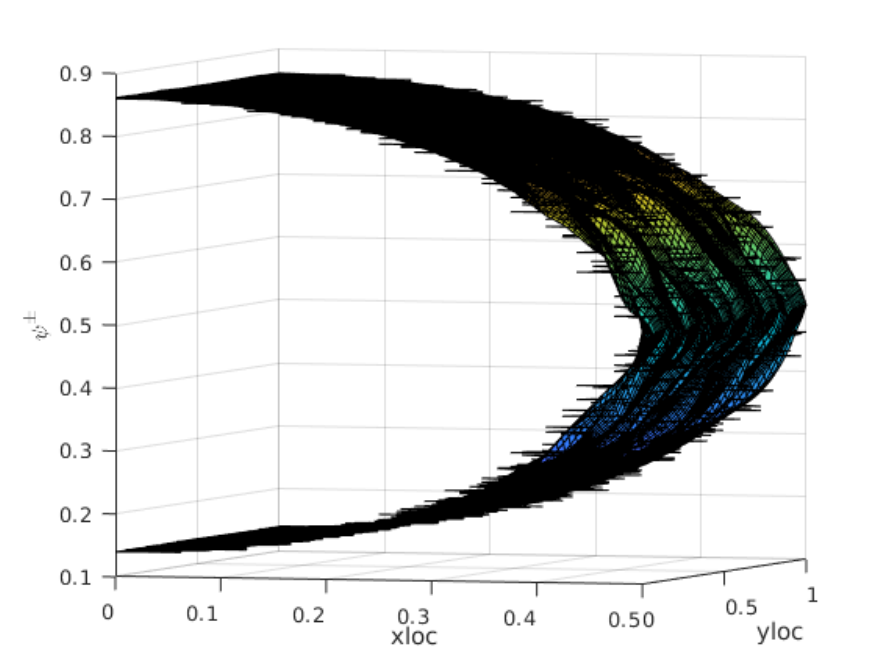}
\caption{\textit{test3}:\textit{Slice plot} of the solution on the plane $\omega^{\pm}$ compared to function $\psi^{\pm}$.}
\label{fig:test3C}
\end{minipage}
\end{figure}

The last numerical test, labelled \textit{test3}, considers a fracture ending inside the domain. The setting is inspired by one of the examples in Section 5.2 of \cite{Angot2009}, re-written as a 3D problem, as before. The domain is a unit edge cube with a vertex in the origin of a reference system $xyz$, Figure~\ref{fig:test3}. The fracture $F$ lies on the plane $x=0.5$ and extends from $0\leq y\leq 0.5$, $0\leq z\leq 0.5$. Homogeneous Dirichelet boundary conditions are prescribed on cube face, edges and vertices lying on the plane $x=0$ and a unit value Dirichelet condition is set on the face, edge and vertices lying on plane $x=1$. Neumann boundary conditions equal to zero are prescribed on the rest of $\partial\D$ and on the whole $\partial F$. Parameters are $K^\D=1$, $K^F=10^{-7}d$, and $\eta=10^{-7}/d$, being $d=10^{-2}$ the thickness of the original equi-dimensional fracture. The solution is obtained on a 3D mesh counting about $4.9\times10^3$ degrees of freedom for $h^\D$ whereas $428$ degrees of freedom are used for $h^F$ and $783$ for $\psi^\star$, $\star={+,-,F}$. For this problem $12$ iterations of the CG solver are required to reach a relative residual of $10^{-7}$. The solution is reported in Figure~\ref{fig:test3} in 3D, whereas Figure~\ref{fig:test3A}, left, shows a \textit{slice plot} on the plane $z=0.5$. A reference solution is computed also for this problem solving a 2D equi-dimensional problem on the plane $z=0.5$, in which the width of the fracture is explicitly represented. The mesh for the reference solution is refined inside fracture width, where maximum element diameter is $10^{-3}$, and is coarser outside, with a maximum diameter of $3\times 10^{-2}$, resulting in about $5\times10^3$ unknowns. The reference solution is shown in Figure~\ref{fig:test3A}, right, and we can see that it is in good agreement with the solution obtained with the proposed method. This is also highlighted by the line plots of Figure~\ref{fig:test3B}, where the two solutions are displayed, overlapped, on six segments placed on the plane $z=0.5$ for different values of $y=\{0.4, 0.5, 0.6, 0.7, 0.8, 0.9\}$. In particular, the computed solution is discontinuous across the fracture, whereas is continuous around $F$, thanks to the chosen enrichment function. Finally, Figure~\ref{fig:test3C} reports the \textit{slice plots} of the 3D solution $h^\D$ on $\omega^\pm$ compared to the corresponding control variables $\psi^\pm$, showing again the expected matching. 

\section{Conclusions}
A method for simulating fractures and barriers in porous media with the DFM model is presented. The method is based on an optimization based domain decomposition strategy, already proposed in previous works, and here modified to account for non permeable fractures and discontinuous solutions in the 3D domain. Key aspects of the method are the use of 3D meshes non conforming to the interfaces and the introduction of a five-field domain decomposition method to decoupled the bulk 3D problem from the problem on the fracture. The eXtended Finite Element Method (XFEM) is adopted to correctly reproduce the discontinuities deriving from filtration-like coupling conditions at the interface between the porous matrix and the fracture. 
The algebraic system deriving from the discretization of the problem with the proposed method is well posed, independently of the choice of the discrete spaces for the involved variables. This is an important aspect in view of its application to complex geometrical configurations.
Four proof-of-concept numerical examples are provided, some of them inspired by reference solutions available in literature and adapted to the present case. The obtained results show the viability of the proposed method and its accuracy with respect to reference solutions obtained with conventional approaches on fine meshes capable of reproducing fracture thickness.

\section*{Acknowledgements}
Partial financial support was provided by INdAM-GNCS. This publication is part of the project NODES which has received funding from the MUR-M4C2 1.5 of PNRR with grant agreement no. ECS00000036

\section*{References}
\bibliographystyle{elsarticle-num}
\bibliography{scico2mp,dfn,misc}

\end{document}